\documentclass[reqno]{amsart}
\usepackage[usenames, dvipsnames]{color}
\usepackage{amsthm}
\usepackage{amsmath}
\usepackage{amssymb}
\usepackage{amscd}
\usepackage{graphics}
\usepackage{latexsym}
\usepackage{stmaryrd}
\usepackage{xspace}
\usepackage[frame,ps,dvips,matrix,arrow,curve,rotate]{xy}
\usepackage{amsfonts}
\usepackage{mathrsfs}

\theoremstyle{plain}
\newtheorem{theorem}{Theorem}[section]
\newtheorem{thm}[theorem]{Theorem}

\newtheorem{lem}[theorem]{Lemma}
\newtheorem{prop}[theorem]{Proposition}
\newtheorem{claim}[theorem]{Claim}

\theoremstyle{definition}

\newtheorem{conj}[theorem]{Conjecture}

\newtheorem{hyp}[theorem]{Hypothesis}

\theoremstyle{remark}

\definecolor{titlecol}{named}{BrickRed}
\definecolor{headcol}{named}{Violet}
\definecolor{seccol}{named}{Red}
\definecolor{sseccol}{named}{Bittersweet}
\definecolor{pbcol}{named}{Black}
\definecolor{sncol}{named}{Brown}
\definecolor{acol1}{named}{Red}
\definecolor{acol2}{named}{Apricot}

\newcommand{\CC}{\mathbb{C}}

\newcommand{\PP}{\mathbb{P}}

\newcommand{\mc}{\mathcal}

\newcommand{\OO}{\mc{O}}

\newcommand{\M}{\overline{\mc{M}}}

\newcommand{\m}{\overline{\text{M}}}

\newcommand{\kgnb}[1]{\m_{#1}}
\newcommand{\Kgnb}[1]{\M_{#1}}

\newcommand{\lt}{\left}
\newcommand{\rt}{\right}

\newsavebox{\sembox}
\newlength{\semwidth}
\newlength{\boxwidth}

\newsavebox{\semrbox}
\newlength{\semrwidth}
\newlength{\boxrwidth}

\def\mC{{\mathbb C}}

\def\w{\bigwedge\nolimits}
\def\O{{\mathcal O}}
\def\Hilb{{\rm Hilb}}
\def\P{\mathbb{P}}

\title{Rational surfaces in Index-one Fano hypersurfaces}

\author[Beheshti]{Roya Beheshti}
\address{Department of Mathematics and Statistics \\ 
Queen's University \\ 
Kingston, Ontario, Canada \\
K7L 3N6}
\email{beheshti@mast.queensu.ca} 

\author[Starr]{Jason Michael Starr}
\address{Department of Mathematics \\
Massachusetts Institute of Technology \\
Cambridge MA 02139}
\email{jstarr@math.mit.edu} 

\date{\today}

\begin{document}


\begin{abstract}
  We give the first evidence for a conjecture that a general,
  index-one, Fano hypersurface is not unirational: (i) a general point
  of the hypersurface is contained in no rational surface ruled,
  roughly, by low-degree rational curves, and (ii) a general point is
  contained in no image of a Del Pezzo surface.
\end{abstract}


\maketitle



\section{Introduction}

For complex, projective varieties a classical notion is
unirationality: A variety rationally dominated by projective space is
\emph{unirational}.  A modern notion is rational connectedness: A
variety is \emph{rationally connected} if every pair of points is
contained in a rational curve.  Every unirational variety is
rationally connected.  The 2 notions agree for curves and surfaces.
Conjecturally they disagree in higher dimensions.

\begin{conj} \label{conj-IM}
For every integer $n\geq 4$ there exists a non-unirational, smooth,
degree-$n$ hypersurface in $\PP^n$.
\end{conj}

A smooth hypersurface in $\PP^n$ of degree $d\leq n$ is an
\emph{index-$(n+1-d)$, Fano manifold}.  By ~\cite{C} and ~\cite{KMM},
every Fano manifold is rationally connected.  Versions of Conjecture
~\ref{conj-IM} have been around for decades.  The specific case $n=4$
is attributed to Iskovskikh and Manin, ~\cite{IskManin}.

In ~\cite{Ksimple}, Koll\'ar suggested an approach to proving
Conjecture ~\ref{conj-IM}.  A general point of an $n$-dimensional,
unirational variety is contained in a $k$-dimensional, rational
subvariety for each $k < n$.  Thus, Conjecture ~\ref{conj-IM} for
$n\geq 5$ follows from the next conjecture (which fails for $n=4$).

\begin{conj} \label{conj-surf}
For every integer $n\geq 5$, there exists a smooth, degree-$n$
hypersurface in $\PP^n$ whose general point is contained in no
rational surface.
\end{conj}  

We give the first evidence for Conjecture ~\ref{conj-surf}.

\begin{thm} \label{thm-1}
For every integer $n\geq 5$, every smooth, degree-$n$ hypersurface $X$
in $\PP^n$ contains a countable union of closed, codimension-$2$
subvarieties containing the image of every generically-finite,
rational transformation $\PP^1 \times \PP^1 \dashrightarrow X$ mapping
a general fiber $\{t\}\times \PP^1$ isomorphically to an
$(n-1)$-normal, smooth, rational curve in $X$.
\end{thm}

\begin{thm} \label{thm-2}
For every integer $n\geq 5$, every smooth, degree-$n$ hypersurface in
$\PP^n$ contains a countable union of closed, codimension-$2$
subvarieties containing the image of every generically finite, regular
morphism from a Del Pezzo surface to $X$.
\end{thm}

A projective variety is \emph{$k$-normal} if every global section of
the restriction of $\OO_{\PP^n}(k)$ is the restriction of a global
section on $\PP^n$.  

There are 2 approaches.  First, given a rational surface $S$ and a
regular morphism $f:S\rightarrow X$, to prove deformations of $f$ are
contained in a codimension-$2$ subvariety of $X$, it suffices to prove
$\bigwedge^{n-4}(f^*T_X/T_S)/\text{Torsion}$ has only the zero
section.  In Section ~\ref{sec-normal}, this is used to prove Theorems
~\ref{thm-1} and ~\ref{thm-2}.

Second, a rational surface with a pencil of rational curves gives a
morphism from $\PP^1$ to a parameter space of rational curves on $X$.
There is a canonical construction of algebraic differential forms on
the parameter space. Since $\PP^1$ has only the zero form, these forms
limit rational curves on the parameter space.  In Section
~\ref{sec-forms}, this is used to prove the following generalization
of Theorem ~\ref{thm-1}.

\medskip\noindent
\begin{thm} \label{thm-3}
For every integer $n\geq 5$, every smooth, degree-$n$ hypersurface $X$
in $\PP^n$ contains a countable union of closed, codimension-$2$
subvarieties containing the image of every generically-finite,
rational transformation $S \dashrightarrow X$ from a surface with a
pencil of curves mapping the general curve isomorphically to an
$(n-1)$-normal, smooth curve of genus $0$ or $1$, also assumed
non-degenerate if the genus is $1$.
\end{thm}

As the second approach does not apply to Theorem ~\ref{thm-2}, the
first approach is more productive.  However, further progress in
proving Conjecture ~\ref{conj-surf} will likely use both approaches,
as well as new ideas.

\bigskip

\paragraph{{\bf Acknowledgments.}} 
The authors thank Mike Roth for many illuminating discussions. 

\bigskip




\section{The first approach} \label{sec-normal}

Let $X$ be a smooth, degree-$n$ hypersurface in $\PP^n$, $n\geq 5$.  
Denote by $\Hilb(X)$ the Hilbert scheme of $X$.
Theorem ~\ref{thm-1} follows easily from the next theorem.

\begin{thm}\label{hmain}
Let $Z$ be an irreducible subscheme of
$\text{Hom}(\PP^1,\Hilb(X))$ satisfying,
\begin{enumerate}
\item[(i)]
the associated morphism $Z\times \PP^1 \rightarrow \Hilb(X)$
does not factor through the projection $Z\times \PP^1 \rightarrow Z$,
and
\item[(ii)]
the image of a general point of $Z\times \PP^1$ parametrizes a smooth,
$(n-1)$-normal, rational curve in $X$.
\end{enumerate}
Then there exists a codimension $\geq 2$ subvariety of $X$ containing
all curves parametrized by $Z\times \PP^1$.
\end{thm}

A morphism $\PP^1 \rightarrow \Hilb(X)$ is equivalent to a
closed subscheme $S'\subset \PP^1 \times X$, flat over $\PP^1$.
If a general point of $\PP^1$ parametrizes a
smooth rational curve, then $S'$ is an irreducible surface.  Any
desingularization $S$ of $S'$ is a surface fitting in a diagram,

\begin{equation} \label{S}
\xymatrix{ S \ar[r]^f \ar[d]^{\pi} & X \\ \P^1
}
\end{equation}

\begin{prop} \label{sprop}
Let $Z$ be an irreducible subvariety of
$\text{Hom}(\PP^1,\Hilb(X))$ satisfying,
\begin{enumerate}
\item[(i)]
the associated morphism $Z\times \PP^1 \rightarrow \Hilb(X)$
does not factor through the projection $Z\times \PP^1 \rightarrow Z$,
\item[(ii)]
a general point of $Z\times \PP^1$ parametrizes a smooth curve in $X$,
and
\item[(iii)]
there is no codimension $2$ subvariety of $X$ containing all curves
parametrized by $Z\times \PP^1$.
\end{enumerate}
Then, for the morphism $\PP^1\rightarrow \Hilb(X)$ parametrized
by a general point of $Z$, the torsion-free sheaf
$\bigwedge^{n-4}(f^*T_X/T_S)/\text{Torsion}$ has a nonzero global
section.
\end{prop}

\begin{proof}
Replacing $Z$ by a dense, Zariski open subset if necessary, assume $Z$
is smooth.  Let $V \subset Z \times \P^1 \times X$ be the pullback of
the universal family to $Z \times \P^1$


Let $g: \widetilde{V} \to V$ be a desingularization of $V$.  Denote by
$\phi$ the projection map from $V$ to $Z \times X$, and denote by $p$
the projection map from $V$ to $Z$.  Let $\widetilde{\phi}= \phi \circ
g$, and let $\tilde{p} = p \circ g$.

$$\xymatrix{\widetilde{V} \ar@(ru, lu)[rr]^{\tilde{g}} 
\ar[r]^{g} \ar[ddr]_{\tilde{p}} & 
V \ar[d] \ar@(dr, ur)[dd]^{p} \ar[r]^-{\phi} & Z \times X \\ 
& Z \times \P^1 \ar[d] \\
& Z. }$$

Replacing $Z$ by a dense, Zariski open subset if necessary, assume
$\tilde{p}$ is smooth, cf. ~\cite[Corollary III.10.7]{H}.  
Associated to the morphism $\tilde{g}$ is the derivative map,
$$
d\tilde{g}: T_{\widetilde{V}} \rightarrow \tilde{g}^* T_X.
$$
Associated to the morphism $\tilde{p}$ is the derivative map,
$$
d\tilde{p}: T_{\widetilde{V}} \rightarrow \tilde{p}^* T_Z.
$$
By hypothesis, $d\tilde{p}$ is surjective.  Denote by $T_{\tilde{p}}$
the kernel of $d\tilde{p}$.  Because $Z\times \PP^1 \rightarrow
\Hilb(X)$ does not factor through $Z$, the restriction of
$\tilde{g}$ to a general fiber of $\tilde{p}$ maps generically
finitely to its image.  Therefore the following map is generically
injective,
$$
d\tilde{g}: T_{\tilde{p}} \rightarrow \tilde{g}^* T_X.
$$
As $T_{\tilde{p}}$ is locally free and $\widetilde{V}$ is integral, in
fact $d\tilde{g}$ is injective.  Denoting by $\mc{N}$ the cokernel of
$\tilde{g}^* T_X$ by $d\tilde{g}(T_{\tilde{p}})$, there is a
commutative diagram with exact rows
$$
\xymatrix{
0 \ar[r] & 
T_{\tilde{p}} \ar[r] \ar[d]_{=} & 
T_{\widetilde{V}} \ar[r]^{d\tilde{p}} \ar[d]^{d\tilde{g}} &
\tilde{p}^*T_Z \ar[r] \ar[d]^{u} & 
0 \\
0 \ar[r] & 
T_{\tilde{p}} \ar[r] &
\tilde{g}^* T_X \ar[r] &
\mc{N} \ar[r] &
0}
$$
By generic smoothness, the rank of $d\tilde{g}$ at a general point of
$\widetilde{V}$ equals the dimension of the closure of
$\text{Image}(\tilde{g})$.  By hypothesis, this is $\geq n-2$.
Therefore the rank of $u$ at a general point is $\geq n-4$.  Thus, the
restriction of $u$ to a general $(n-4)$-plane in the fiber of
$\tilde{p}^* T_Z$ has rank $n-4$.  A general $(n-4)$-plane is the
tangent space of a general $(n-4)$-dimensional subvariety of $Z$.
Therefore, after replacing $Z$ by the smooth locus of a general
$(n-4)$-dimensional subvariety of $Z$, assume $Z$ is
$(n-4)$-dimensional and $u$ is generically injective.  

Associated to $u$, there is an induced map,
$$
\bigwedge^{n-4}(u): \tilde{p}^* \bigwedge^{n-4} T_Z \rightarrow
\bigwedge^{n-4} \mc{N}/\text{Torsion}.
$$
Because $u$ is generically injective and $n\geq 5$, this map is
generically injective.  

Let $z$ be a general point of $Z$, and denote by $S$ and
$\widetilde{S}$ the fibers of $p$ and $\tilde{p}$ over $z$,
respectively.  Since $\widetilde{V}$ is smooth, $\widetilde{S}$ is a
smooth surface.  Let $f: S \to X$ be the restriction of $\phi$ to $S$,
and let $\tilde{f}: \widetilde{S} \to X$ be the corresponding map from
$\widetilde{S}$.  The restriction of $\mc{N}$ to $S$ is precisely $f^*
T_X/T_S$.  The restriction of $\tilde{p}^* T_Z$ to $S$ is precisely
the trivial vector bundle $T_{Z,z}\otimes_\CC \OO_S$.  Since $z$ is
general, the restriction of $\bigwedge^{n-4}(u)$ is generically
injective.  Therefore, it is
a nonzero map,
$$
\bigwedge^{n-4}(u)|_S: (\bigwedge^{n-4} T_{Z,z})\otimes_\CC \OO_S \rightarrow
(\bigwedge^{n-4} f^* T_X/T_S)/\text{Torsion}.
$$
Since $T_{Z,z}$ is $(n-4)$-dimensional, this is equivalent to a
nonzero global section of $(\bigwedge^{n-4} f^*
T_X/T_S)/\text{Torsion}$ (well-defined up to nonzero scaling).
\end{proof}

\begin{prop} \label{sprop2}
Let $\PP^1 \rightarrow \Hilb(X)$ be a morphism with associated
diagram as in Equation~\ref{S}.  If the curve parametrized by a general
point of $\PP^1$ is smooth and $(n-1)$-normal then
$$
h^0(S,\omega_S\otimes \bigwedge^{n-2}f^* T_X) = 0.
$$
\end{prop}

\begin{proof}
Pulling back the short exact sequence of tangent bundles
$$ 0 \to T_X \to T_{\P^n}|_X \to N_{X/\P^n} \cong \O_X(n) \to 0$$
to $S$ and taking its $(n-1)^{\text st}$ exterior power gives 
another short exact sequence 
$$ 0 \to \bigwedge^{n-1} f^*T_X \to \w^{n-1} f^* T_{\P^n} \to 
f^*\O_X(n) \otimes \w^{n-2} f^*T_X \to 0.$$
Tensoring this sequence with $\omega_S \otimes 
f^*\O_X(-n)$ gives the following short exact sequence
$$
0 \to \omega_S \otimes f^*\O_X(-n) \otimes 
\bigwedge^{n-1} f^*T_X   \to 
\omega_S \otimes f^*\O_X(-n) \otimes 
\bigwedge^{n-1} f^* T_{\P^n} \to 
\omega_S \otimes \bigwedge^{n-2} f^*T_X \to 0.
$$

Applying the long exact sequence of cohomology, 
$h^0(S, \omega_S \otimes \w^{n-2} f^*T_X)$ equals $0$ if both, 
\begin{enumerate}
\item[(i)]
$h^0(S, \omega_S \otimes f^*\O_X(-n) \otimes 
\w^{n-1} f^*T_{\P^n})$ equals $0$, and
\item[(ii)]
$h^1(S,  \omega_S \otimes f^*\O_X(-n) \otimes 
\w^{n-1} f^*T_X)$ equals $0$.
\end{enumerate}

\textbf{Proof of (i).} 
Consider the Euler exact sequence on $\P^n$
$$ 0 \to \O_{\P^n} \to \O_{\P^n}(1)^{n+1}\to T_{\P^n} \to 0.$$
Pulling this back to $S$, and taking its $n^{\text{th}}$ 
exterior power gives the following exact sequence
$$ 0 \to \w^{n-1} f^* T_{\P^n}   
\to f^*\O_X(n)^{\oplus(n+1)} \to f^*\O_X(n+1) \to 0.$$
Tensoring with $\omega_S \otimes 
f^*\O_X(-n)$ gives an injective map 
$$
\omega_S \otimes f^*\O_X(-n) \otimes 
\w^{n-1} f^*T_{\P^n} \to 
\omega_S^{\oplus(n+1)}.
$$
Thus it suffices to prove $h^0(S,\omega_S)$ equals $0$, which follows
from the hypothesis that $S$ is a rational surface.

\textbf{Proof of (ii).}
There is a canonical isomorphism
$$ \omega_S \otimes f^*\O_X(-n) \otimes \w^{n-1} f^*T_X \cong \omega_S
\otimes f^*\O_X(-n+1).$$  
So by Serre duality, it suffices to prove
$h^1(S, f^*\O_X(n-1))$ equals $0$. Let $C$ be a general fiber of the map 
$\pi: S \to \P^1$. There is a short exact sequence
\begin{equation}\label{eq0}
0 \to f^*\O_X(n-1) \otimes I_{C/S} \to 
f^* \O_X(n-1) \to f^* \O_X(n-1) |_C \to 0,
\end{equation}
where $\mathcal I_{C/S}$ is the ideal sheaf of $C$ in $S$. 
By hypothesis, the image of $C$ by $f$ is 
$(n-1)$-normal in $\P^n$, therefore the map 
$$
H^0(S, f^*\O_X(n-1)) 
\to H^0(C, f^*\O_X(n-1)|_C)
$$  
is surjective.  The long exact sequence 
of cohomology to the sequence in Equation \ref{eq0} gives an isomorphism
\begin{equation} \label{eqhalf}
H^1(S, f^*\O_X(n-1) \otimes I_{C/S}) \cong 
H^1(S, f^*\O_X(n-1)).
\end{equation}

Because $S$ is a smooth surface and the general fiber of $\pi$ is a
smooth, rational curve, $R^1\pi_*\mc{F}$ is the zero sheaf for every
$\pi$-relatively globally-generated, coherent $\OO_S$-module $\mc{F}$.
Because $f^*\O_{X}(n-1)$ is globally-generated, 
it is $\pi$-relatively globally-generated.  Because
$I_{C/S}$ equals $\pi^* \O_{\P^1}(-1)$, the twist $f^*\O_X(n-1)\otimes
I_{C/S}$ is $\pi$-relatively globally-generated.
Thus $R^1\pi_*(f^*\O_X(n-1))$ and $R^1\pi_*(f^*\O_X(n-1)\otimes
I_{C/S})$ are each the zero sheaf.  So,
by the Leray spectral sequence, there are canonical isomorphisms
\begin{equation}\label{eq1}
H^1(S, f^*\O_X(n-1)) \cong 
H^1(\P^1, \pi_*(f^*\O_X(n-1))),
\end{equation}
\begin{equation}\label{eq2}
H^1(S, f^*\O_X(n-1) \otimes \mathcal I_{C/S}) \cong 
H^1(\P^1, \pi_*(f^*\O_X(n-1) \otimes \mathcal I_{C/S})).
\end{equation}

Taken together, Equations \ref{eqhalf}, \ref{eq1} and \ref{eq2} give a
canonical isomorphism 
$$ H^1(\P^1, \pi_*(f^*\O_X(n-1))) \cong 
H^1(\P^1, \pi_*(f^*\O_X(n-1)) \otimes \O_{\P^1}(-1)).$$
This is possible only if $h^1(\P^1, \pi_*(f^*\O_X(n-1)))$ equals $0$.
\end{proof}

\begin{proof}[Proof of Theorem \ref{hmain}] 
Let $Z$ satisfy the hypotheses of Proposition \ref{sprop}, and let $S$
and $f$ satisfy the conclusion of Proposition \ref{sprop}.  The
injective map $df:T_S \rightarrow f^* T_X$ induces a multiplication
map,
$$
\bigwedge^2 T_S \otimes \bigwedge^{n-4} f^* T_X \rightarrow
\bigwedge^{n-2} f^* T_X.
$$
The image sheaf is precisely $\w^2 T_X \otimes (\w^{n-4}(f^*
T_X/T_S))/\text{Torsion}$.  Tensoring with the canonical bundle of
$\omega_S$, this gives an injective map
$$ (\w^{n-4}(f^*T_X/T_S))/\text{Torsion} 
\to \omega_S \otimes \w^{n-2} f^*T_X.$$
By hypothesis, $(\w^{n-4}(f^*T_X/T_S))/\text{Torsion}$ has a nonzero
global section.  Therefore $\omega_S \otimes \w^{n-2} f^*T_X$ also has
a nonzero global section.

On the other hand, for $Z$ satisfying the hypothesis of
Theorem~\ref{hmain}, Proposition \ref{sprop2} implies,
$$
h^0(S,\omega_S \otimes \w^{n-2} f^*T_X) = 0.
$$
Thus $Z$ does not satisfy the hypothesis of Proposition \ref{sprop2},
i.e., it does not satisfy Hypothesis (iii).  Therefore there exists a
codimension $2$ subvariety of $X$ containing all the curves
parametrized by $Z\times \PP^1$.
\end{proof}

\begin{proof}[Proof of Theorem~\ref{thm-1}]
For every generically-finite, rational transformation $\PP^1 \times
\PP^1 \dashrightarrow X$ restricting to a closed immersion on a
general fiber, there is an associated rational transformation
$$
\PP^1 \dashrightarrow \Hilb(X), \ \ t \mapsto
\text{Image}(\{t\}\times \PP^1).
$$
By properness of the Hilbert scheme and the valuative criterion, this
extends to a regular morphism.  Therefore, associated to each rational
transformation is an element in the Hom-scheme
$\text{Hom}(\PP^1,\Hilb(X))$.  Those rational transformations
satisfying the hypothesis of Theorem ~\ref{thm-1} give a locally
closed subset of $\text{Hom}(\PP^1,\Hilb(X))$.  As
$\text{Hom}(\PP^1,\Hilb(X))$ is a countable union of quasi-projective
varieties, this subset is also a countable union of quasi-projective
subvarieties.  By Theorem ~\ref{hmain}, for each such subvariety $Z$,
there is a codimension-$2$ subvariety of $X$ containing every curve
parametrized by $Z\times \PP^1$.  This subvariety contains the image
of each rational transformation $\PP^1 \times \PP^1 \dashrightarrow X$
giving a point in $Z$.  Therefore, there exists a countable union of
codimension-$2$ subvarieties of $X$ containing the image of every
rational transformation satisfying the hypothesis of
Theorem~\ref{hmain}.
\end{proof}

The proof of the Theorem \ref{thm-2} is similar to the proof of
Theorem \ref{hmain}.  There is a preliminary proposition.

\begin{prop} \label{DP}
Let $X$ be a smooth hypersurface of degree $n$ in $\P^n$.  For every
Del Pezzo surface $S$ and every generically finite morphism $f: S \to
X$, the only global section of
$\w^{n-4} (f^*T_X/T_S)/\text{Torsion}$ is the zero section.
\end{prop}

\begin{proof}
The proof is similar to the proof of Theorem \ref{hmain}. 
By the same type of argument, it suffices to prove,
\begin{enumerate}
\item[(i)]
$h^0(S, \omega_S \otimes f^*\O_X(-n) \otimes 
\w^{n-1} f^*T_{\P^n})$ equals $0$, and
\item[(ii)]
$h^1(S,  \omega_S \otimes f^*\O_X(-n) \otimes 
\w^{n-1} f^*T_X)$ equals $0$.
\end{enumerate}

The proof of (i) is the same as in the proof of Theorem \ref{hmain},
since $h^0(S, \omega_S)$ equals $0$.  As for (ii), there is a
canonical isomorphism
$$\omega_S \otimes f^*\O_X(-n) \otimes 
\w^{n-1}f^*T_X \cong (\omega_S^{-1} \otimes 
f^*\O_X(n-1))^{-1}.$$
Denote $\omega_S^{-1}\otimes f^*\O_X(n-1)$ by $L$.  The sheaf
$f^*\O_X(n-1)$ is globally generated.  By the hypothesis that $S$ is a
Del Pezzo surface, $\omega_S^{-1}$ is ample.  Thus $L$ is ample.  By
Kodaira vanishing, $h^1(S,L^{-1})$ equals $0$.  So, using the
canonical isomorphism, $h^1(S, \omega_S \otimes f^*\O_X(-n) \otimes
\w^{n-1}f^*T_X)$ equals $0$.
\end{proof}

\begin{proof}[Proof of Theorem \ref{thm-2}]
By the same countability argument as at the beginning of the section,
it suffices to prove that for every flat family
$$ \xymatrix{ D \ar[d]^p \ar[r]^-{\phi} & X \times V \\
V }$$
such that $\phi$ is generically finite and a general fiber of $p$ is a
Del Pezzo surface, the image of $\phi$ is contained in a subvariety of
codimension $\geq 2$.  Let $S$ be the fiber of $p$ over a general
point of $V$, and let $f$ be the restriction $\phi|_S: S \to X$.  As
in the proof of the Proposition \ref{sprop}, if the image of $\phi$ is
contained in no subvariety of codimension $\geq 2$, then $H^0(S,
\w^{n-4} (f^*T_X/T_S)/\text{Torsion})$ has a nonzero global section.
Thus Proposition \ref{DP} proves the image of $\phi$ is contained in a
subvariety of codimension $\geq 2$.
\end{proof}


\section{The second approach} \label{sec-forms}

Let $X$ be a smooth hypersurface in $\PP^n$ of degree $n$. We 
denote by $\Kgnb{g,n}(X)$ the Kontsevich moduli stack of families 
of genus-$g$, $n$-pointed, stable maps to $X$.  
The associated coarse moduli space is denoted by $\kgnb{g,n}(X)$,  
cf. \cite{FP}. 

For every integral subscheme $M$ of $\kgnb{g,0}(X)$, denote by
$X(M)$ the smallest closed subvariety of $X$ containing every curve
parametrized by $M$.

\noindent 
\begin{hyp} \label{hyp-XM}
Let $M$ be an integral, closed subscheme of $\kgnb{g,0}(X)$.
\begin{enumerate}
\item[(i)]
The curves parametrized by $M$ are contained in no codimension $2$
subscheme of $X$, i.e., $\text{dim}(X(M)) \geq n-2$.
\item[(ii)]
The integer $g$ equals $0$ or $1$.
\item[(iii)]
A general point of $M$ parametrizes an embedded, smooth,
$(n-1)$-normal curve.  If $g$ equals $1$, also the curve is
nondegenerate.  
\item[(iv)]
The dimension of $M$ equals $\text{dim}(X)-2 = n-3$.
\end{enumerate}
\end{hyp}

\medskip

\begin{thm} \label{thm-nonuni1}
For every integer $n\geq 4$, for every smooth, degree-$n$ hypersurface $X$
in $\PP^n$, and for every integral closed subscheme $M$ of
$\kgnb{g,0}(X)$ satisfying Hypothesis \ref{hyp-XM},  
every desingularization of $M$ has a nonzero canonical form.
In particular, $M$ is not uniruled.
\end{thm}

There are 2 components of the proof: a global construction of certain
$(n-3)$-forms on $M$ following ~\cite{dJS2}, and a local description
of the forms proving they are nonzero.  The construction is
~\ref{lem-nonuni3}, the local description is Lemma~\ref{lem-nonuni5},
and the nonvanishing result is Claim~\ref{claim-nonuni1}.

Let $\widetilde{M}$ be a finite type scheme and let $\nu:\widetilde{M}
\rightarrow \kgnb{g,0}(X)$ be a morphism.  Later, 
$\widetilde{M}$ will be a desingularization of an integral
subscheme $M$ of $\kgnb{g,0}(X)$.
For every integer $p$, ~\cite[Corollary 4.3]{dJS2} gives a map,
$$
\psi_p:H^1(X,\Omega^{p+1}_X) \rightarrow
H^0(\Kgnb{g,0}(X),\Omega^p_{\Kgnb{g,0}(X)}). 
$$
Denote,
$$
\widetilde{M}_{\text{stack}} := \lt(
\widetilde{M}\times_{\kgnb{g,0}(X)}\Kgnb{g,0}(X)
\rt)_{\text{red}},
$$
i.e., the associated reduced stack of the 2-fibered product.  There
are projections,
$$
\pi_1:\widetilde{M}_{\text{stack}} \rightarrow \widetilde{M},
$$
$$
\pi_2:\widetilde{M}_{\text{stack}} \rightarrow \Kgnb{g,0}(X).
$$
There are associated pullback maps on $p$-forms,
$$
\pi_1^*:H^0(\widetilde{M},\Omega^p_{\widetilde{M}}) \rightarrow
H^0(\widetilde{M}_{\text{stack}},\Omega^p_{\widetilde{M}_{\text{stack}}}),
$$
$$
\pi_2^*:H^0(\Kgnb{g,0}(X),\Omega^p_{\Kgnb{g,0}(X)})
\rightarrow
H^0(\widetilde{M}_{\text{stack}},\Omega^p_{\widetilde{M}_{\text{stack}}}).
$$
Hence there is a map,
$$
\pi_2^*\circ \psi:H^1(X,\Omega^{p+1}_X) \rightarrow
H^0(\widetilde{M}_{\text{stack}},\Omega^{p}_{\widetilde{M}_{\text{stack}}}).
$$
Does this map factor through $\pi_1^*$, i.e., does there exist a map,
$$
\phi_p: H^1(X,\Omega^{p+1}_X) \rightarrow
H^0(\widetilde{M},\Omega^{p}_{\widetilde{M}}),
$$
such that $\pi_1^*\circ \phi_p$ equals $\pi_2^*\circ \psi_p$?

\begin{lem} \label{lem-nonuni3}
If $\widetilde{M}$ is smooth there is a unique map of $\mC$-vector spaces,
$$
\phi_p:
H^1(X,\Omega^{p+1}_X) \rightarrow
H^0(\widetilde{M},\Omega^{p}_{\widetilde{M}}),
$$
such that $\pi_1^*\circ \phi_p$ equals $\pi_2^*\circ \psi_p$.
\end{lem}

\begin{proof}
Since $\widetilde{M}$ is smooth and $\widetilde{M}_{\text{stack}}$ is
reduced, 
\cite[Proposition 3.6]{dJS2} implies the pullback map,
$$
H^0(\widetilde{M},\Omega^{p}_{\widetilde{M}}) \rightarrow
H^0(\widetilde{M}_{\text{stack}},
\Omega^{p}_{\widetilde{M}_{\text{stack}}}),
$$
is an isomorphism.  
\end{proof}

Let $C$ be a smooth curve in $X$ with corresponding point $[C]$ in
$\Kgnb{g,0}(X)$.  For every integer $p$,
restriction to the fiber at $[C]$ defines a map,
$$
\psi_{p,[C]}:H^1(X,\Omega^{p+1}_X) \rightarrow
\Omega^{p}_{\Kgnb{g,0}(X)}|_{[C]}.
$$
The Zariski tangent space to $\Kgnb{g,0}(X)$ at
$[C]$ is $H^0(C,N_{C/X})$.  The dual vector space is the fiber of
$\Omega_{\Kgnb{g,0}(X)}$ at $[C]$.  The $p^{\text{th}}$ exterior
power is the fiber of $\Omega^{p}_{\Kgnb{g,0}(X)}$ at $[C]$.  
Therefore $\psi_{p,[C]}$ is equivalent to a linear map,
$$
\psi_{p,[C]}:H^1(X,\Omega^{p+1}_X) \rightarrow
\text{Hom}(\bigwedge^{p}H^0(C,N_{C/X}), \mC).
$$
What is the map $\psi_{p,[C]}$?  In other words, 
for an element in $H^1(X,\Omega^{p+1}_X)$, what is the associated
$p$-linear alternating map on $H^0(C,N_{C/X})$?

When $X$ is a smooth hypersurface in $\PP^n$ and $p=n-3$, 
the answer follows as in ~\cite[Theorem
  5.1]{dJS2}.  For a smooth hypersurface $X$ in $\PP^n$,
Griffiths computed the cohomology groups $H^q(X,\Omega^p_X)$, 
cf. ~\cite[Section 8]{GR}.  In particular, there is an exact sequence,
\begin{equation} \label{eqn-Gr}
H^0(X,\Omega^{n-1}_{\PP^n}(X)|_X)\rightarrow
H^0(X,\Omega^n_{\PP^n}(2X)|_X) \rightarrow H^1(X,\Omega^{n-2}_X)
\rightarrow 0.
\end{equation}
Therefore every element of $H^1(X,\Omega^{n-2}_X)$ is the image
$\overline{\beta}$ of an element $\beta$ in
$H^0(X,\Omega^n_{\PP^n}(2X)|_X)$.  

There is a short exact sequence of locally free
$\OO_C$-modules,
$$
0 \rightarrow N_{C/X} \rightarrow N_{C/\PP^n} \rightarrow
N_{X/\PP^n}|_C \rightarrow 0. 
$$
Taking the $(n-2)^{\text{nd}}$ exterior power gives a short exact
sequence,
$$
0 \rightarrow \bigwedge^{n-2} N_{C/X} \rightarrow \bigwedge^{n-2}
N_{C/\PP^n} \rightarrow 
(\bigwedge^{n-3} N_{C/X}) \otimes N_{X/\PP^n}|_C \rightarrow 0.
$$
Twisting each term by $N_{X/\PP^n}^\vee|_C$ gives an exact sequence,
\begin{equation} \label{eqn-ses}
0 \rightarrow \bigwedge^{n-2} N_{C/X} \otimes N_{X/\PP^n}^\vee|_C \rightarrow
\bigwedge^{n-2} N_{C/\PP^n} \otimes N_{X/\PP^n}^\vee|_C \rightarrow
\bigwedge^{n-3} N_{C/X} \rightarrow 0.
\end{equation}
Applying the long exact sequence of cohomology, there is a connecting
map,
$$
H^0(C,\bigwedge^{n-3} N_{C/X}) \rightarrow H^1(C,\bigwedge^{n-2}
N_{C/X} \otimes N_{X/\PP^n}^\vee|_C).
$$

Now $\bigwedge^{n-2}N_{C/X}$ is the determinant of $N_{C/X}$, which is
canonically isomorphic to $\omega_C\otimes (\Omega^{n-1}_X)^\vee|_C$.
By adjunction, $\Omega^{n-1}_X$ is isomorphic to
$\Omega^n_{\PP^n}(X)|_X$.  Also, $N_{X/\PP^n}$ is isomorphic to
$\OO_{\PP^n}(X)|_X$.  Putting this together gives a canonical
isomorphism,
$$
\bigwedge^{n-2} N_{C/X}\otimes N_{X/\PP^n}^\vee|_C \cong \omega_C
\otimes (\Omega^n_{\PP^n}(2X)^\vee)|_C.
$$
Serre duality gives an isomorphism,
$$
H^1(C,\omega_C \otimes (\Omega^n_{\PP^n}(2X)^\vee)|_C) \cong
H^0(C,\Omega^n_{\PP^n}(2X)|_C)^\vee.
$$
The pullback map $H^0(X,\Omega^n_{\PP^n}(2X)|_X) \rightarrow
H^0(C,\Omega^n_{\PP^n}(2X)|_C)$ determines a transpose map,
$$
H^0(C,\Omega^n_{\PP^n}(2X)|_C)^\vee \rightarrow
H^0(X,\Omega^n_{\PP^n}(2X)|_X)^\vee.
$$
Finally, every element $\beta$ of $H^0(X,\Omega^n_{\PP^n}(2X)|_X)$
determines a linear functional,
$$
H^0(X,\Omega^n_{\PP^n}(2X)|_X)^\vee \rightarrow \mC.
$$
Putting all this together, every element $\beta$ of
$H^0(X,\Omega^n_{\PP^n}(2X)|_X)$ determines a linear functional,
$$
\widetilde{\beta}:H^0(C,\bigwedge^{n-3} N_{C/X}) \rightarrow \mC.
$$

\begin{lem} \label{lem-nonuni5}
Let $X$ be a smooth hypersurface in $\PP^n$.  For every element
$\beta$ of $H^0(X,\Omega^n_{\PP^n}(2X)|_X)$,
$\phi_{n-3,[C]}(\overline{\beta})$ equals the restriction of
$\widetilde{\beta}$ to $\bigwedge^{n-2} H^0(C,N_{C/X})$, up to nonzero
scaling.
\end{lem}

\begin{proof}
This follows by a diagram-chase.  Here are the main points.
There is a commutative diagram with exact rows and columns,
\begin{equation} \label{eqn-CD}
\xymatrix{
& 0 \ar[d] & 0 \ar[d] & \\
& T_C \ar[d] \ar[r]^{\text{Id}} & T_C \ar[d] \\
0 \ar[r] & T_X|_C \ar[d] \ar[r] & T_{\PP^n}|_C \ar[d] \ar[r] & 
N_{X/\PP^n}|_C \ar[d]^{\text{Id}} \ar[r] & 0 \\
0 \ar[r] &  N_{C/X} \ar[d] \ar[r] & N_{C/\PP^n} \ar[d] \ar[r] & 
N_{X/\PP^n}|_C \ar[r] & 0 \\ 
&  0 &  0 &  
}
\end{equation}

The following three invertible sheaves are isomorphic,
$$
\omega_C \otimes \bigwedge^{n-2} N_{C/X}^\vee \cong \Omega^{n-1}_X|_C
\cong \Omega^n_{\PP^n}(X)|_C.
$$
Denote any by $L$.  Twisting
Equation~\ref{eqn-CD} by $L$ gives a commutative diagram with exact
rows and columns,
{\small
\begin{equation} \label{eqn-CD2}
\xymatrix@C=14pt{
& 0 \ar[d] & 0 \ar[d] \\
0 \ar[r] & \bigwedge^{n-2} N_{C/X}^\vee \ar[d] \ar[r]^-{\cong} & 
\bigwedge^{n-1} N_{C/\PP^n}^\vee \otimes \OO_{\PP^n}(X)|_C \ar[d] 
\ar[r] &  0 \ar[d] \\
0 \ar[r] & \Omega^{n-2}_X|_C \ar[d] \ar[r] & 
\Omega^{n-1}_{\PP^n}(X)|_C \ar[d] \ar[r] & 
\Omega^{n-1}_X|_C\otimes \OO_{\PP^n}(X)|_C \ar[d]^{\cong} 
\ar[r] & 0 \\
0 \ar[r] & \omega_C \otimes \bigwedge^{n-3} N_{C/X}^\vee \ar[d] 
\ar[r] & \omega_C
\otimes \bigwedge^{n-2} N_{C/\PP^n}^\vee \otimes \OO_{\PP^n}(X)|_C 
\ar[d] \ar[r] & 
\omega_C \otimes \bigwedge^{n-2} N_{C/X}^\vee \otimes
\OO_{\PP^n}(X)|_C 
\ar[d] \ar[r] & 0 \\ 
& 0 & 0 & 0
}
\end{equation}
} 
For a sheaf $\mc{E}$ on $X$, $\mc{E}(X)|_X$ denotes the tensor product
$\mc{E}\otimes \OO_{\PP^n}(X)|_X$.  And for a sheaf $\mc{F}$ on $C$,
$\mc{F}(X)|_C$ denotes the tensor product, $\mc{F}\otimes
\OO_{\PP^n}(X)|_C$.

\medskip\noindent
Consider the last map in the first column of Equation~\ref{eqn-CD2}.
There is an associated map of cohomology groups $H^1(C,-)$,
$$
H^1(C,\Omega_X^{n-2}|_C)\rightarrow H^1(C,\omega_C\otimes
(\bigwedge^{n-3} N_{C/X})^\vee).
$$
By Serre duality this is equivalent to,
$$
H^1(C,\Omega_X^{n-2}|_C)\rightarrow
\text{Hom}(H^0(C,\bigwedge^{n-3}N_{C/X}),\mC).
$$

There is a natural multiplication map,
$$
\bigwedge^{n-3}H^0(C,N_{C/X})\rightarrow
H^0(C,\bigwedge^{n-3}N_{C/X}),
$$
having transpose,
$$
\text{Hom}(H^0(C,\bigwedge^{n-3}N_{C/X}),\mC) \rightarrow
\text{Hom}(\bigwedge^{n-3} H^0(C,N_{C/X}),\mC).
$$
Composing the previous map with the transpose gives a map,
$$
H^1(C,\Omega_X^{n-2}|_C)\rightarrow \text{Hom}(\bigwedge^{n-3}
H^0(C,N_{C/X}),\mC).
$$
Composing the restriction map with the last map gives a map,
$$
H^1(X,\Omega^{n-2}_X) \rightarrow
\text{Hom}(\bigwedge^{n-3} H^0(C,N_{C/X}),\mC).
$$
As in the proof of ~\cite[Theorem 5.1]{dJS2}, this equals the map defined
in Lemma ~\ref{lem-nonuni3}.  

On the other hand, the second row of Equation~\ref{eqn-CD2}
is the short exact sequence giving rise to Equation~\ref{eqn-Gr}.
Therefore, the associated map,
$$
H^0(X,\Omega^n_{\PP^n}(2X)|_X) \rightarrow \text{Hom}(\bigwedge^{n-3}
H^0(C,N_{C/X},\mC),
$$
defining $\phi_{n-3,[C]}(\overline{\beta})$ is obtained by taking the
push-out of the second row of Equation~\ref{eqn-CD2} by the last map
in the first column.  This push-out is canonically isomorphic to the
third row of Equation~\ref{eqn-CD2}.  Therefore
$\phi_{n-3,[C]}(\overline{\beta})$ comes from the connecting map in
cohomology associated to the third row of Equation~\ref{eqn-CD2}

Finally, the third row of Equation~\ref{eqn-CD2} is obtained from the
exact sequence in Equation~\ref{eqn-ses} by dualizing and tensoring
with $\omega_C$.  In particular, using Serre duality, the connecting
map for the third row,
$$
H^0(C,\omega_C\otimes \bigwedge^{n-2} N_{C/X}^\vee \otimes
\OO_{\PP^n}(X)|_C) \rightarrow H^1(C,\omega_C \otimes \bigwedge^{n-3}
N_{C/X}^\vee),
$$
equals the transpose of the map,
$$
H^0(C,\bigwedge^{n-3} N_{C/X}) \rightarrow H^1(C,\bigwedge^{n-2}
N_{C/X} \otimes N_{X/\PP^n}^\vee|_C).
$$
Since this is the map used to define $\widetilde{\beta}$, the
restriction of $\widetilde{\beta}$ equals
$\phi_{n-3,[C]}(\overline{\beta})$.
\end{proof}

\begin{proof}[Proof of Theorem~\ref{thm-nonuni1}]
Let $M$ be an integral, closed subscheme of $\kgnb{g,0}(X)$
satisfying Hypotheses \ref{hyp-XM}.  Let
$\nu:\widetilde{M}\rightarrow M$ be a desingularization of $M$.  There
is an open dense subscheme $V$ of $\widetilde{M}$ on which $\nu$ is
unramified and over which the pullback family of stable maps is a
family of smooth, embedded, $(n-1)$-normal curves, also assumed
nondegenerate if $g$ equals $1$.  For every point $[C]$ in $V$,
$T_{V,[C]}$ is a subspace of $H^0(C,N_{C/X})$.  Therefore
$\bigwedge^{n-3} T_{V,[C]}$ is a $1$-dimensional subspace of
$\bigwedge^{n-3} H^0(C,N_{C/X})$.  By Hypothesis \ref{hyp-XM}(i), for
a general point $[C]$ of $V$, the image of $\bigwedge^{n-3} T_{V,[C]}$
in $H^0(C,\bigwedge^{n-3} N_{C/X})$ is nonzero.

Associated to the short exact sequence in Equation~\ref{eqn-ses},
there is an exact sequence of cohomology,
\begin{equation} \label{eqn-ses2}
H^0(C,\bigwedge^{n-2} N_{C/\PP^n}\otimes N_{X/\PP^n}^\vee|_C)
\rightarrow H^0(C,\bigwedge^{n-3} N_{C/X}) \xrightarrow{\delta}
H^1(C,\bigwedge^{n-2} N_{C/X}\otimes N_{X/\PP^n}^\vee|_C).
\end{equation}

\begin{claim} \label{claim-nonuni1}
$h^0(C,\bigwedge^{n-2} N_{C/\PP^n}\otimes
N_{X/\PP^n}^\vee|_C)$ equals $0$.
\end{claim}

To prove Claim~\ref{claim-nonuni1}, 
first observe that $\bigwedge^{n-2} N_{C/\PP^n}$
is isomorphic to $N_{C/\PP^n}^\vee \otimes \bigwedge^{n-1}
N_{C/\PP^n}$.  By adjunction, $\bigwedge^{n-1} N_{C/\PP^n}$ is
isomorphic to $\omega_C \otimes (\Omega^n_{\PP^n})^\vee|_C \cong
\omega_C \otimes \OO_{\PP^n}(n+1)|_C$.  Since $N_{X/\PP^n}^\vee$ is
isomorphic to $\OO_{\PP^N}(-n)|_X$, altogether there is an
isomorphism,
$$
N_{C/\PP^n}\otimes N_{X/\PP^n}^\vee|_C \cong \omega_C\otimes
N_{C/\PP^n}^\vee\otimes 
\OO_{\PP^n}(1)|_C. 
$$
Now $N_{C/\PP^n}^\vee$ is a subsheaf of $\Omega_{\PP^n}|_C$.
Therefore $H^0(C,\omega_C \otimes N_{C/\PP^n}^\vee \otimes
\OO_{\PP^n}(1)|_C)$ is a subspace of $H^0(C,\omega_C \otimes
\Omega_{\PP^n}(1)|_C)$.  There is an exact sequence,
$$
0 \to \Omega_{\PP^n}(1) \to H^0(\PP^n,\OO_{\PP^n}(1))\otimes_k
\OO_{\PP^n} \to \OO_{\PP^n}(1) \to 0,
$$
where the last map is the canonical map.  By Hypothesis
~\ref{hyp-XM}(ii), $g$ equals $0$ or $1$.  If $g$ equals $0$,
$h^0(C,\omega_C)$ equals $0$, and thus $h^0(C,\omega_C\otimes
\Omega_{\PP^n}(1)|_C)$ also equals $0$.  If $g$ equals $1$, then
$\omega_C$ is isomorphic to $\OO_C$.  Therefore $H^0(C,\omega_C
\otimes \Omega_{\PP^n}(1)|_C)$ is the kernel of the following map,
$$
H^0(\PP^n,\OO_{\PP^n}(1)) \rightarrow H^0(C,\OO_{\PP^n}(1)|_C).
$$
By Hypothesis \ref{hyp-XM}(iii), $C$ is nondegenerate, and thus the
kernel is trivial.  This proves Claim~\ref{claim-nonuni1}, both when
$g$ equals $0$ and when $g$ equals $1$.

Because of Claim~\ref{claim-nonuni1}, the map $\delta$ from
Equation~\ref{eqn-ses2} is injective.  As in the proof of
Lemma~\ref{lem-nonuni5}, the target of $\delta$ is canonically
isomorphic to
$H^0(C,\Omega_{\PP^n}(2X)|_C)^\vee$.
Now $\Omega_{\PP^n}(2X)$ is isomorphic to $\OO_{\PP^n}(-n-1)\otimes
\OO_{\PP^n}(2n) \cong \OO_{\PP^n}(n-1)$.  By Hypothesis
~\ref{hyp-XM}(ii), $C$ is $(n-1)$-normal.  Thus the following map is
injective
$$
H^0(C,\Omega_{\PP^n}(2X)|_C)^\vee \rightarrow
H^0(X,\Omega_{\PP^n}(2X)|_X)^\vee.
$$
In particular, there exists an element $\beta$ of
$H^0(X,\Omega_{\PP^n}(2X)|_X)$ whose associated map
$\widetilde{\beta}$ is nonzero on $\bigwedge^{n-3} T_{V,[C]}$.  By
Lemma~\ref{lem-nonuni5}, the $(n-3)$-form $\phi_p(\overline{\beta})$
is nonzero at $[C]$.  Therefore $\phi_p(\overline{\beta})$ is a
nonzero, canonical form on $\widetilde{M}$, proving $\widetilde{M}$ is
not uniruled, cf.~\cite[Corollary 4.12]{De}.
\end{proof}

\begin{proof}[Proof of Theorem~\ref{thm-3}]
By a proof similar to that of Theorem ~\ref{thm-1}, it suffices to
prove there exists no uniruled subvariety $M$ of $\kgnb{g,0}(X)$
satisfying Hypothesis~\ref{hyp-XM}(i)--(iii).

Let $M$ be an integral, uniruled subvariety of $\kgnb{g,0}(X)$
satisfying Hypothesis~\ref{hyp-XM}(i) and (ii), and whose general point
parametrizes a smooth, embedded curve.  
Since $M$ is uniruled, there exists a quasi-projective variety $Z$
and a dominant, generically-finite morphism
$$
Z\times \PP^1 \rightarrow M.
$$
As in the proof of Proposition ~\ref{sprop}, for a
general $(n-4)$-dimensional subvariety $Z'$ of $Z$, 
the closed image of $Z'\times \PP^1$
is an $(n-3)$-dimensional, uniruled subvariety $M'$ of $M$ satisfying
Hypothesis~\ref{hyp-XM}(i) and (ii) and containing a general point of
$M$.
So $M'$ is uniruled and satisfies Hypothesis~\ref{hyp-XM} (i), (ii)
and (iv).  By Theorem~\ref{thm-nonuni1}, $M'$ does not satisfy (iii).
Since $M'$ contains a general point of $M$, also $M$ does not satisfy
(iii).  
Therefore every integral, closed subvariety $M$ satisfying
Hypothesis~\ref{hyp-XM}(i)--(iii) is not uniruled.
\end{proof}



\newcommand{\closer}{\vspace{-1.5ex}}

\end{document}